\documentclass[a4paper,10pt,openany]{article}
\usepackage{indentfirst}
\usepackage[greek,french,english]{babel}
\usepackage[iso-8859-7]{inputenc}

\usepackage{color,listings,titlesec}
\usepackage[headings]{fullpage}
\usepackage{amsfonts}
\usepackage{fancyhdr}
\usepackage{graphicx}
\usepackage{amsmath}
\newtheorem{Proposition}{Proposition}

\newtheorem{Definition}[Proposition]{Definition}

\newtheorem{Theorem}[Proposition]{Theorem}

\newcommand{\g}{\mathfrak{g}}
\newcommand{\ad}{ad}
\newcommand{\h}{\mathfrak{h}}

\begin{document}
\begin{center}

{\Large \textbf{Deformation quantization and invariant differential operators.}}
\vspace{1.5cm}

\textbf{Panagiotis Batakidis}\footnote{Department of Mathematics, Aristotle University of Thessaloniki.}

\end{center}

\textbf{Abstract.}  In this note we explain how the techniques of deformation quantization in the sense of Kontsevich can be used to describe the algebra of invariant differential operators on Lie groups.

\textbf{Version fran\c caise abr\'eg\'ee.}

Soit $\g$ une alg\`ebre de Lie, $\h\subset\g$ une sous-alg\`ebre et $\lambda$ un  caract\`ere de $\h$. En prenant comme motivation la quantification par d\'eformation de Kontsevich et la g\'en\'eralization de Cattaneo-Felder on regarde $\h_{\lambda}^{\bot}$ comme sous-vari\'et\'e coisotrope de la vari\'et\'e de Poisson $\g^{\ast}$.  C'est l'approche de Cattaneo-Torossian aux alg\`ebres de Lie pour calculer dans  l'alg\`ebre d\'eform\'ee $\left(U_{(\epsilon)}(\g)/U_{(\epsilon)}(\g)\h_{\lambda}\right)^{\h}$. Nous \'etudions \'egalement motiv\'es par les Conjectures de Duflo et Corwin-Greenleaf. Le d\'epart est l'isomorphisme $(U(\mathfrak{g})/U(\mathfrak{g})\mathfrak{h}_{\lambda})^{\mathfrak{h}}\simeq \mathbb{D}(\mathfrak{g},\mathfrak{h},\lambda)$ de Koornwider qui interpr\`ete l'alg\`ebre  $(U(\mathfrak{g})/U(\mathfrak{g})\mathfrak{h}_{\lambda})^{\mathfrak{h}}$ comme les op\'erateurs lin\'eaires qui laissent invariant l'espace $C^{\infty}(G,H,\lambda)$ des fonctions complexes  $\theta : G\longrightarrow \mathbb{C}$ qui satisfaient $\theta (g\cdot \exp X)=e^{-i\lambda(X)}\theta(g), \;\;\forall X\in \mathfrak{h}, \forall g\in G$.

Nous  annon\c cons qu'l ya un isomorphisme non-canonique $\overline{\beta}_{\mathfrak{q},(\epsilon)}\circ\partial_{q_{(\epsilon)}^{\frac{1}{2}}}\circ \overline{T}_1^{-1}T_2:\; H^0_{(\epsilon)}(\h^{\bot}_{\lambda},d^{(\epsilon)}_{\h^{\bot}_{\lambda},\mathfrak{q}})\stackrel{\simeq}{\longrightarrow} \left(U_{(\epsilon)}(\g)/U_{(\epsilon)}(\g)\h_{\lambda}\right)^{\h}$ entre l'alg\`ebre de r\'eduction sur l'espace affine $\h_{\lambda}^{\bot} :=\{f\in\mathfrak{g}^{\ast}/f|_{\mathfrak{h}}=-\lambda\}$ et l'alg\`ebre d\'eform\'ee $\left(U_{(\epsilon)}(\g)/U_{(\epsilon)}(\g)\h_{\lambda}\right)^{\h}$. Ensuite on \'etudie la sp\'ecialization $H^0_{(\epsilon=1)}(\h^{\bot}_{\lambda},d^{(\epsilon=1)}_{\h^{\bot}_{\lambda},\mathfrak{q}})$, l'alg\`ebre de r\'eduction $H^0(\h_{\lambda}^{\bot},d_{\h_{\lambda}^{\bot},\mathfrak{q}})$ d\'efinie sans le param\`etre de d\'eformation $\epsilon$, et  $H^0(\h_{t\lambda}^{\bot},d_{\h_{t\lambda}^{\bot},\mathfrak{q}}), t\in\mathbb{R}$ l'alg\`ebre de r\'eduction sans $\epsilon$ et en d\'eformant le  caract\`ere $\lambda$. 

Nous comparons les objets correspondants de la part de $H^0_{(\epsilon)}(\h^{\bot}_{\lambda},d^{(\epsilon)}_{\h^{\bot}_{\lambda},\mathfrak{q}})$ avec les objets associ\'es du cot\'e $\left(U_{(\epsilon)}(\g)/U_{(\epsilon)}(\g)\h_{\lambda}\right)^{\h}$ \c ca veut dire $\mathbb{D}_{(T=1)}(\g,\h,\lambda)$ et $\left((U(\g)/U(\g)\h_{t\lambda})^{\h}\right)$. On arrive finalement d'annon\c cer que

 \[\mathbb{D}_{(T=1)}(\g,\h,\lambda)\stackrel{alg}{\simeq}H^0_{(\epsilon=1)}(\h_{\lambda}^{\bot},d^{(\epsilon=1)}_{\h_{\lambda}^{\bot},\mathfrak{q}}).\]

\section{Introduction.}
Let $G$ be a nilpotent, connected and simply connected Lie group, $H\subset G$ a subgroup and $\lambda\in\g^{\ast}$ a character of $\h$. Let $C^{\infty}(G,H,\lambda)$ be the vector space of $C^{\infty}$ functions $\theta : G\longrightarrow \mathbb{C}$ that satisfy the condition
\[ \theta (g\cdot \exp X)=e^{-i\lambda(X)}\theta(g), \;\;\forall X\in \mathfrak{h}, \forall g\in G\]
and $\mathbb{D}(\mathfrak{g},\mathfrak{h},\lambda)$ the algebra of linear differential operators that leave the space $C^{\infty}(G,H,\lambda)$ invariant and commute with the left translation on $G$,  
\[D(C^{\infty}(G,H,\lambda))\subset C^{\infty}(G,H,\lambda),\;\;\text{and}\;\;D(L(g)\theta)=L(g)(D(\theta)).\]
Let also $m(\tau)$ be the multiplicities in the spectral decomposition of the representation $\tau_{\lambda}$
\[\tau_{\lambda} \cong \int_{\widehat{G}} m(\tau)\tau d\mu(\tau)\cong\int_{(f+\mathfrak{h}^{\bot})/H}\tau_l d\nu(l)\]
where $\tau_{\lambda}$ is the representation $Ind(G,H,\lambda)=L^2(G,H,\lambda)$ and $\widehat{G}$ is the space of irreducible representations, or otherwise the unitary dual of $G$. Furthermore $\nu$ is a finite and positive measure equivalent to the Lebesgue measure on the affine space $\lambda+\h^{\bot}$. Finally $\mu=K_{\ast}(\nu)$ (where $K:\g^{\ast}\longrightarrow \widehat{G}$ is the Kirillov map) is our measure on $\widehat{G}$.

The \textbf{\textlatin{Corwin-Greenleaf}} conjecture says that
\[\mathbb{D}(\mathfrak{g},\mathfrak{h},\lambda)\; \text{is commutative} \;\;  \mathbf{ \Leftrightarrow}\;\; m(\tau)<+\infty\;\; \mu- \text{a.e}.\]
 To state a second conjecture of our interest, let $\mathfrak{g}$ be a nilpotent Lie algebra, $\mathfrak{h}\subset\mathfrak{g}$ a subalgebra and $\lambda$ a character of $\mathfrak{h}$. Let $U(\mathfrak{g})\mathfrak{h}_{\lambda}$ be the left ideal of $U(\mathfrak{g})$ generated by the family $<X+\lambda(X),\;X\in\mathfrak{h}>$.
 
Then the \textbf{\textlatin{Duflo}} conjecture says that
\[\text{The algebras}\; C_{poiss}((S(\mathfrak{g})/S(\mathfrak{g})\mathfrak{h}_{\lambda})^{\mathfrak{h}})\; \text{and}\; C_{ass}[(U(\mathfrak{g})/U(\mathfrak{g})\mathfrak{h}_{\lambda})^{\mathfrak{h}}]\;\;\text{are isomorphic.}\]
The first center refers to the Poisson algebra structure of $(S(\mathfrak{g})/S(\mathfrak{g})\mathfrak{h}_{\lambda})^{\mathfrak{h}}$ while the second to the associative structure of $(U(\mathfrak{g})/U(\mathfrak{g})\mathfrak{h}_{\lambda})^{\mathfrak{h}}$.

The analytic interest of the second conjecture stems from a theorem of \textlatin{Koornwider} stating that $(U(\mathfrak{g})/U(\mathfrak{g})\mathfrak{h}_{\lambda})^{\mathfrak{h}}\simeq \mathbb{D}(\mathfrak{g},\mathfrak{h},\lambda)$. Taking into consideration the interpretation of $U(\g)$ as the linear differential operators on $G$, the Duflo conjecture is actually a question about the structure of the invariant differential operators on the homogeneous space $G/H$.

In this note we will prove that $(U(\mathfrak{g})/U(\mathfrak{g})\mathfrak{h}_{\lambda})^{\mathfrak{h}}\simeq \mathbb{D}(\mathfrak{g},\mathfrak{h},\lambda)$ is isomorphic to the reduction algebra related to these data, an algebra playing a central role to the deformation quantization theory. we will also provide results of the same nature for related reduction algebras. The results of this note were presented in \cite{BAT}. 

\section{Deformation quantization and generalizations.}
We start by reminding the standard result
\begin{Theorem}[\cite{K}]\label{Kontsevich}
Let $\pi$ be a Poisson bivector of $\mathbb{R}^k$ and $F,G\in C^{\infty}(\mathbb{R}^k)$. 
The product
\[F \ast_KG:=F\cdot G+\sum_{n=1}^{\infty}\epsilon^n\left(\frac{1}{n!}\sum_{\Gamma\in\mathbf{Q_{n,\bar{2}}}}\omega_{\Gamma}B_{\Gamma,\pi}(F,G)\right)\]
is associative.
\end{Theorem}
The set $\mathbf{Q_{n,\bar{2}}}$ is a special family of graphs $\Gamma$. To every admissible graph $\Gamma$ corresponds a bidifferential operator $B_{\Gamma}(F,G):=\sum_{R,S}b_{\Gamma}^{R,S}\partial_R(F)\partial_S(G)$ on $C^{\infty}(\mathbb{R}^k)\times C^{\infty}(\mathbb{R}^k)$. Note that the functions $b_\Gamma^{R,S}$ depend on $\Gamma$ and are $n$-linear in the bivector $\pi$. The coefficient $\omega_{\Gamma}\in\mathbb{R}$ is computed by integrating a differential form $\Omega_{\Gamma}$ (which is also encoded in $\Gamma$) over a concentration manifold 

\[\hat{C}^+_{n,\bar{m}}=\{(z_1,\dots,z_n,z_{\bar{1}},\dots,z_{\bar{m}}) / z_i\in\mathbb C, \mathfrak{Im}(z_i)>0, z_{\bar{i}}\in\mathbb R, z_{\bar{i}}<z_{\bar{j}}\; \text{for}\; i<j\}/G_2,\]
where $G_2$ is the 2-dimensional Lie group of dilations $\langle z_k\mapsto az_k+b, a>0, b\in\mathbb R\rangle$. So $\hat{C}^+_{n,\bar{m}}\subset (H^+)^n\times R^m$.

In the general setting that we briefly recalled, the central result of \cite{K} is the following
\begin{Theorem}[\cite{K}]\label{Kontsevich 2}
Let $\mathcal{U}:\;\mathcal{T}_{poly}(\mathbb{R}^k)\longrightarrow \mathcal{D}_{poly}(\mathbb{R}^k)$ be the map defined by the Taylor coefficients:
\[\mathcal{U}_n:=\sum_{\overline{m}\geq 0}\left(\sum_{\Gamma\in\mathbf{Q}_{n,\overline{m}}}\omega_{\Gamma}B_{\Gamma}\right).\]
Then $\mathcal{U}$ is an $L_{\infty}-$ morphism and a quasi-isomorphism between the differential graded Lie algebras (DGLA) $\mathcal{T}_{poly}(\mathbb{R}^k)$ of polyvector fields on $\mathbb{R}^k$ and $\mathcal{D}_{poly}(\mathbb{R}^k)$ of polydifferential operators on $\mathbb{R}^k$.
\end{Theorem}
The key point in the proof is a Stokes equation integrating the form $\Omega_{\Gamma}$ on $\hat{C}_{n,\overline{m}}^+$. This way Kontsevich also reached the Duflo isomorphism (applying these to the case of the Poisson manifold $\g^{\ast}$). 
Our approach is using the following generalization of what is so far said.
Let $X$ be a Poisson manifold, and $C\subset X$ a coisotropic submanifold. Let $NC$ be the normal bundle, $\mathcal{A}=\bigoplus_{i=0}^{rank(\pi)}\mathcal{A}^i$ with $\mathcal{A}^i=\Gamma(C,\wedge^iNC)$ be the graded commutative algebra of sections of the exterior algebra of the normal bundle $NC$ and $\tilde{\mathcal{D}}(\mathcal{A})=\oplus_n\tilde{\mathcal{D}}^n(\mathcal{A})$ with $\tilde{\mathcal{D}}^n(\mathcal{A})=\prod_{p+q-1=n}Hom^p(\otimes^q\mathcal{A},\mathcal{A})$. Finally we set  \[\mathcal{B}=\bigoplus_{j=0}^{\infty}\mathcal{B}^j,\;\; \mathcal{B}^0=\Gamma(C,S((NC)^{\ast})), \;\mathcal{B}^j=0,\;\textlatin{if}\;j\neq 0.\]

\begin{Theorem}[\cite{CF3}]
Consider the DGLA $\mathcal{T}(X,C)$ and $\tilde{\mathcal{D}}(\mathcal{A})$ of polyvector fields and polydifferential operators on $X$. There is an $L_{\infty}-$ quasi-isomorphism
\[\mathcal{F}_R:\;\mathcal{T}(X,C)\longrightarrow \tilde{\mathcal{D}}(\mathcal{A}),\]
whose first Taylor coefficient $\mathcal{F}_R^1$ is the composition
\[F_{HKR}\circ \hat{F}:\;\mathcal{T}(X,C)\simeq \tilde{\mathcal{T}}(\mathcal{B})\;\stackrel{\hat{F}}{\rightarrow}\; \tilde{\mathcal{T}}(\mathcal{A})\;\stackrel{F_{HKR}}{\rightarrow}\;\tilde{\mathcal{D}}(\mathcal{A}).\]
\end{Theorem}
where $F_{HKR}$ is the Hochschild-Konstant-Rosenberg map
\[F_{HKR}:= S_{\mathcal{A}}(Der(\mathcal{A})[-1])\longrightarrow \oplus_{j=0}^{\infty}Hom_{\mathbf{K}}(\otimes^j\mathcal{A},\mathcal{A})\] and $\hat{F}$ is a kind of Fourier transform (for more details we refer to \cite{CF3}).

We will now apply the previous results for $\g^{\ast}$ with the natural Poisson structure coming from the Lie structure and the orthogonal space $\h^{\bot}$ of a subalgebra $\h\subset\g$ as a coisotropic submanifold. Some modifications in the theory are needed since the conjectures of the introduction refer to the affine space $-\lambda+\h^{\bot}\subset \g^{\ast}$.

\section{Applications for non-commutative harmonic analysis.}

\textbf{Reduction equations.}
Let $\mathfrak{h}^{\bot}:=\{l\in\mathfrak{g}^{\ast}/ l(\mathfrak{h})=0\}$, $\mathfrak{h}_{\lambda}^{\bot}:=\{f\in\mathfrak{g}^{\ast}/f|_{\mathfrak{h}}=-\lambda\}$ and $\mathfrak{q}$ a supplementary space of $\h$ in $\g$. Consider the differential $d^{(\epsilon)}_{\mathfrak{h}^{\bot},\mathfrak{q}}:\;S(\mathfrak{q})[\epsilon]\longrightarrow S(\mathfrak{q})[\epsilon]\otimes\mathfrak{h}^{\ast}$ where $d^{(\epsilon)}_{\mathfrak{h}^{\bot},\mathfrak{q}}:=\sum_{i=1}^{\infty}\epsilon^i d^{(i)}_{\mathfrak{h}^{\bot},\mathfrak{q}}$ and $d^{(i)}_{\mathfrak{h}^{\bot},\mathfrak{q}}=\sum_{\Gamma\in \mathcal{B}_i\cup \mathcal{BW}_i}\overline{\omega}_{\Gamma}B_{\Gamma}$.

The second sum is on two families of graphs in $Q_{1,\bar{i}}$: the set ${\mathcal B}_i$ of Bernoulli graphs and the set ${\mathcal {BW}}_i$ of Bernoulli-attached-to-a-wheel graph (see figure 1).
The first terms of the differential $d_{\mathfrak{h}^{\bot},\mathfrak{q}}(F)=0$ are
\[d^1_{\mathfrak{h}^{\bot},\mathfrak{q}}(F_n)=0,\;\;\;\;\; d^3_{\mathfrak{h}^{\bot},\mathfrak{q}}(F_n)+d^1_{\mathfrak{h}^{\bot},\mathfrak{q}}(F_{n-2})=0,\]
\[d^5_{\mathfrak{h}^{\bot},\mathfrak{q}}(F_n)+d^3_{\mathfrak{h}^{\bot},\mathfrak{q}}(F_{n-2})+d^1_{\mathfrak{h}^{\bot},\mathfrak{q}}(F_{n-4})=0\;\ldots\]
For example, $d^1_{\mathfrak{h}^{\bot},\mathfrak{q}}(F_n)=0\Rightarrow F_n\in S(\mathfrak{q})^{\mathfrak{h}}$.

\begin{figure}[h!]
\begin{center}
\includegraphics[width=9cm]{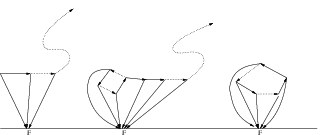}
\caption\footnotesize{A  Bernoulli graph ($\mathcal{B}_3$) in $d_{\mathfrak{h}^{\bot},\mathfrak{q}}^3$, a Bernoulli-attached-to-a-wheel ($\mathcal{BW}_7$) in $d_{\mathfrak{h}^{\bot},\mathfrak{q}}^7$, and a wheel-type graph.}
\end{center}
\end{figure}

\begin{Definition}[\cite{BAT}, $\mathcal{x}$ 2.3](Reduction algebra)
We define the reduction space $H^0_{(\epsilon)}(\h^{\bot},d^{(\epsilon)}_{\h^{\bot},\mathfrak{q}})$ of polynomials in the formal variable $\epsilon$ as the vector space of polynomials - solutions of the system of linear differential equations
\[d^{(\epsilon)}_{\h^{\bot},\mathfrak{q}}(F_{(\epsilon)})=0,\;F_{(\epsilon)}\in S(\mathfrak{q})[\epsilon].\]
The reduction algebra is $H^0_{(\epsilon)}(\h^{\bot},d^{(\epsilon)}_{\h^{\bot},\mathfrak{q}})$ equipped with the product $\ast_{CF,\epsilon}$.
\end{Definition}
We now explain our main results. Let $G$ be a real Lie group, $\h\subset\g$ a subalgebra, $\lambda\in \mathfrak{h}^{\ast}$ a character of $\h$ and $f\in\mathfrak{g}^{\ast}$ s.t $f|_{\mathfrak{h}}=\lambda$. Let also $\mathfrak{q}$ be a supplementary space of $\h$ in $\g$. 

We set the deformed tensor algebra to be $T_{(\epsilon)}(\g):=\mathbf{R}[\epsilon]\otimes T(\g)$. Let $\mathcal{I}_{\epsilon}$ be the two-sided ideal $<X\otimes Y-Y\otimes X-\epsilon[X,Y]>$ of $T_{(\epsilon)}(\g)$. Since our basic model of $U(\g)$ is $T(\g)$ factored by the non-homogeneous ideal $<X\otimes Y-Y\otimes X-[X,Y]>$, we define the deformed universal enveloping algebra of $\g$ as $U_{(\epsilon)}(\g):=T_{(\epsilon)}(\g)/\mathcal{I}_{\epsilon}$.  Let also $\h_{\lambda}:=\{H+\lambda(h),\;H\in\h\}$.

\begin{Theorem}[\cite{BAT}, $\mathcal{x}$ 3.4.3, Theorem 3.1]
The map
\[\overline{\beta}_{\mathfrak{q},(\epsilon)}\circ\partial_{q_{(\epsilon)}^{\frac{1}{2}}}\circ \overline{T}_1^{-1}T_2:\; H^0_{(\epsilon)}(\h^{\bot}_{\lambda},d^{(\epsilon)}_{\h^{\bot}_{\lambda},\mathfrak{q}})\stackrel{\simeq}{\longrightarrow} \left(U_{(\epsilon)}(\g)/U_{(\epsilon)}(\g)\h_{\lambda}\right)^{\h}\]
is an algebra isomorphism.
\end{Theorem}
Here $\overline{\beta}_{\mathfrak{q},(\epsilon)}:\;S(\mathfrak{q})[\epsilon]\longrightarrow U_{(\epsilon)}(\mathfrak{g})/U_{(\epsilon)}(\mathfrak{g})\mathfrak{h}_{\lambda}$ is the symmetrization map, $q(Y) = \det_{\mathfrak{g}} \left(\frac{sinh\frac{\ad Y}{2}}{\frac{\ad Y}{2}}\right)$ and $T_1,T_2$ are differential operators that can be described in terms of Kontsevich graphs \cite{BAT} $\mathcal{x}$ 2.5.2. The theorem is powerful since no condition is needed for the original Lie group $G$.

\textbf{Proof.} We use the $H^0_{(\epsilon)}(\h_{\lambda}^{\bot},d^{(\epsilon)}_{\g^{\ast},\h_{\lambda}^{\bot},\mathfrak{q}})-$ bimodule structure
\[T_1:\;H^0_{(\epsilon)}(\g^{\ast},d^{(\epsilon)}_{\g^{\ast}})\longrightarrow H^0_{(\epsilon)}(\h_{\lambda}^{\bot},d^{(\epsilon)}_{\g^{\ast},\h_{\lambda}^{\bot},\mathfrak{q}}),\;\; G\mapsto G\ast_1 1\] and
\[T_2:\;H^0_{(\epsilon)}(\h^{\bot}_{\lambda},d^{(\epsilon)}_{\h^{\bot}_{\lambda},\mathfrak{q}}) \longrightarrow H^0_{(\epsilon)}(\h_{\lambda}^{\bot},d^{(\epsilon)}_{\g^{\ast},\h_{\lambda}^{\bot},\mathfrak{q}}),\;\;F\mapsto 1\ast_2 F\]
 in the biquantization diagram (see \cite{CKTB} and \cite{CT}) of $\mathfrak{g}^{\ast}$ and $\mathfrak{h}_{\lambda}^{\bot}$.
We show first that $H^0_{(\epsilon)}(\h^{\bot}_{\lambda},d^{(\epsilon)}_{\h^{\bot}_{\lambda},\mathfrak{q}})\subset (U_{(\epsilon)}(\g)/U_{(\epsilon)}(\g)\h_{\lambda})^{\h}$ exploiting results from \cite{CF2},\cite{CF3},\cite{CT} and the bimodule structure we mentioned to pass from  $H^0_{(\epsilon)}(\h^{\bot}_{\lambda},d^{(\epsilon)}_{\h^{\bot}_{\lambda},\mathfrak{q}})$ to $(U_{(\epsilon)}(\g)/U_{(\epsilon)}(\g)\h_{\lambda})^{\h}.$

 For the inverse let $H\in \h$ be a function of first degree $F\in H^0_{(\epsilon)}(\h^{\bot}_{\lambda},d^{(\epsilon)}_{\h^{\bot}_{\lambda},\mathfrak{q}})$. We check the graphs in the expression
 \[(H+\lambda(H))\ast_1 1 \ast_2 F.\]
 
We show that $(H+\lambda(H))\ast_1 (1 \ast_2 F)=[(H+\lambda(H))\ast_1 1] \ast_2 F=0$. 

\begin{figure}[h!]
\begin{center}
\includegraphics[width=9cm]{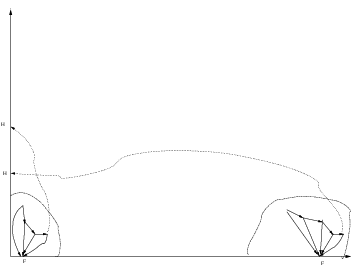}
\caption{\footnotesize Behaviour of the expression $(H+\lambda(H))\ast_1 1 \ast_2 F$ when $s\rightarrow 0$ and $s\rightarrow \infty$.}
\end{center}
\end{figure}

 Then examining the possible and admissible graphs in the concentration manifolds in this expression we end up to the equation 
 \[\sum_{\Gamma^{\alpha}_{\textlatin{int}},\Gamma^{\alpha}_{\textlatin{ext}}}\int_0^{\infty}\hat{\omega}_{\Gamma^{\alpha}_{\textlatin{ext}}}(\textlatin{s})\mathcal{B}_{\Gamma^{\alpha}_{\textlatin{ext}}}\left(\omega_{\Gamma^{\alpha}_{\textlatin{int}}}\mathcal{B}_{\Gamma^{\alpha}_{\textlatin{int}}}\right)\mathrm{d}\textlatin{s}=0\]
 which is equivalent to
 \[\sum_{\alpha}\left(\sum_{\Gamma^{\alpha}_{\textlatin{int}},\Gamma^{\alpha}_{\textlatin{ext}}}\left(\sum_{\textlatin{l},\textlatin{k},\textlatin{m}}\left(B^\textlatin{m}_{\Gamma^{\alpha}_{\textlatin{ext}}}(B^\textlatin{k}_{\Gamma^{\alpha}_{\textlatin{int}}}(\textlatin{F}_\textlatin{l}))\right)\epsilon^{\textlatin{m+k+l}}\right)\right)=0.\]
 This last equation gives (after determining the admissible graphs above) the reduction equations defining $\mathbf{H^0_{(\epsilon)}(\h^{\bot}_{\lambda},d^{(\epsilon)}_{\h^{\bot}_{\lambda},\mathfrak{q}})}$.

\textbf{The specialization algebra} $\mathbf{H^0_{(\epsilon=1)}(\h_{\lambda}^{\bot},d^{(\epsilon=1)}_{\h_{\lambda}^{\bot},\mathfrak{q}})}$.

\begin{Definition}
 The specialization algebra for the affine space $-\lambda+\h^{\bot}=\h_{\lambda}^{\bot}$, is defined as
\[H^0_{(\epsilon=1)}(\h_{\lambda}^{\bot},d^{(\epsilon=1)}_{\h_{\lambda}^{\bot},\mathfrak{q}}):=\left(H^0_{(\epsilon)}(\h_{\lambda}^{\bot},d^{(\epsilon)}_{\h_{\lambda}^{\bot},\mathfrak{q}})/<\epsilon-1>\right).\]
The Cattaneo-Felder product on $H^0_{(\epsilon=1)}(\h_{\lambda}^{\bot},d^{(\epsilon=1)}_{\h_{\lambda}^{\bot},\mathfrak{q}})$ will be also denoted as $\ast_{CF,(\epsilon=1)}$.
\end{Definition}

We have to note here that one may consider also the reduction algebras $H^0(\h_{\lambda}^{\bot},d_{\h_{\lambda}^{\bot},\mathfrak{q}})$ (that is defined without the formal variable $\epsilon$) and $H^0(\h_{t\lambda}^{\bot},d_{\h_{t\lambda}^{\bot},\mathfrak{q}}), t\in\mathbb{R}$ (deforming the character $\lambda$).

Let now $\textlatin{F}^{'}\in  H^0_{(\epsilon)}(\h_{\lambda}^{\bot},d^{(\epsilon)}_{\h^{\bot}_{\lambda},\mathfrak{q}})$. Let $J:\;H^0_{(\epsilon)}(\h_{\lambda}^{\bot},d^{(\epsilon)}_{\h^{\bot}_{\lambda},\mathfrak{q}})\longrightarrow H^0(\h_{\lambda}^{\bot},d_{\h^{\bot}_{\lambda},\mathfrak{q}})$, $J(F^{'}):=\sum_k F^{'}_k$. To describe the image of the map $J$ we have the following:

\begin{Theorem}[\cite{BAT}, $\mathcal{x}$ 3.5.2]
$\bullet$ Let $F\in H^0(\h_{\lambda}^{\bot},d_{\h_{\lambda}^{\bot},\mathfrak{q}})$. Suppose an element $F_{(t)}=\sum_pt^pF_p$ with $F_{(t)}\in H^0(\h_{t\lambda}^{\bot},d_{\h_{t\lambda}^{\bot},\mathfrak{q}})$, $\forall t\in\mathbb{R}^{\ast}$ and $F_{(t=1)}=F$. Let $F_p=\sum_iF_p^{(i)}$ be a decomposition to homogeneous components and $F_{(\epsilon)}:=\epsilon^N\sum F_p^{(i)}\frac{1}{\epsilon^{i+p}}$ ($N>>max(i+p)$). Then $F_{(\epsilon)}\in H^0_{(\epsilon)}(\h_{\lambda}^{\bot},d^{(\epsilon)}_{\h_{\lambda}^{\bot},\mathfrak{q}})$ and $J(F_{(\epsilon)})=F$.

$\bullet$ Let $F\in H^0(\h_{\lambda}^{\bot},d_{\h_{\lambda}^{\bot},\mathfrak{q}})$. Suppose an element $F_{(\epsilon)}=\sum_{0 \geq k\geq n}\epsilon^kF_k$ with $F_{(\epsilon)}\in H^0_{(\epsilon)}(\h_{\lambda}^{\bot},d^{(\epsilon)}_{\h_{\lambda}^{\bot},\mathfrak{q}})$ and $J(F_{(\epsilon)})=F$. Let $F_k=\sum_iF_k^{(i)}$ be a decomposition to homogeneous components and $F_{(t)}:=t^N\sum_{i,k}\frac{1}{t^{i+k}}F_k^{(i)}$ ($N>>max(i+k)$). Then  $\forall t\in\mathbb{R}^{\ast},\;\;F_{(t)}\in H^0(\h_{t\lambda}^{\bot},d_{\h_{t\lambda}^{\bot},\mathfrak{q}})$.
\end{Theorem}

\textbf{Deformations.}  In this section we will make clear the relation between the various reduction algebras presented in the previous part. Using Theorem 5 we will associate them to the appropriate algebras of operators. Let $\g_T:=\g\oplus\mathbb{R}T,\;\;[\g,T]=0$. We note as $\mathcal{P}_{(t)}\left((U(\g)/U(\g)\h_{t\lambda})^{\h}\right)$ the algebra of polynomial families in $t$, $t\longrightarrow u_t\in\left(U(\g)/U(\g)\h_{t\lambda}\right)^{\h}$.
Let $e_t:\;U(\g_T)\longrightarrow U(\g)$ be the surjective map $T\mapsto t$ and $\forall X\in\g,\; X\mapsto X$. Then $\left(U(\g_T)/<T-t>\right)\simeq U(\g)$ and evaluating at $T=t$, we take the surjective map
\[e_{(T=t)}:\;\left( (U(\g_T)/U(\g_T)\h^T_{\lambda})^{\h_T}/<T-t>\right)\hookrightarrow (U(\g)/U(\g)\h_{t\lambda})^{\h}.\]

It turns out that if $t\longrightarrow u_t\in (U(\g)/U(\g)\h_{t\lambda})^{\h}$ is a polynomial family in $t$, then there is a $u_T\in (U(\g_T)/U(\g_T)\h^T_{\lambda})^{\h_T}$ s.t $e_t(u_T)=u_t$.

Now let's examine the notion of specialization from the differential operator point of view.
Let $\mathbb{D}_{(T=1)}(\g,\h,\lambda):= (U(\g_T)/U(\g_T)\h^T_{\lambda})^{\h_T}/<T-~1>$. In the case $t=1$ of the previous example we get
\[\mathbb{D}_{(T=1)}(\g,\h,\lambda)\hookrightarrow (U(\g)/U(\g)\h_{\lambda})^{\h}.\]

More specifically if $u\in \mathbb{D}_{(T=1)}(\g,\h,\lambda)$ then there is a $u_T\in U(\g_T)$ s.t $u=\pi_{(T=1)}(u_T)$. The element $u_t:=e_{(T=t)}(u_T)\in (U(\g)/U(\g)\h_{t\lambda})^{\h}$ defines a polynomial family in $t$ so $u_t\in \mathcal{P}_{(t)}\left((U(\g)/U(\g)\h_{t\lambda})^{\h}\right)$.
Then
\[\mathbb{D}_{(T=1)}(\g,\h,\lambda)\simeq \mathcal{P}_{(t=1)}\left((U(\g)/U(\g)\h_{\lambda})^{\h}\right).\]

We get similar results for specialization from the reduction algebra point of view. 

\begin{Theorem}[\cite{BAT}, $\mathcal{x}$ 3.5.3]
$\bullet$ Let $t\mapsto F_t \in H^0(\h_{t\lambda}^{\bot},d_{\h_{t\lambda}^{\bot},\mathfrak{q}})$ be a polynomial family in $t$. Then there is a $F_T\in  H^0((\h_{\lambda}^T)^{\bot},d_{(\h_{\lambda}^T)^{\bot},\mathfrak{q}})$ s.t $e_t(F_T)=F_t$.

$\bullet$ $\mathbb{D}_{(T=1)}(\g,\h,\lambda)\stackrel{alg}{\simeq}H^0_{(\epsilon=1)}(\h_{\lambda}^{\bot},d^{(\epsilon=1)}_{\h_{\lambda}^{\bot},\mathfrak{q}})$.
\end{Theorem}

\textbf{Acknowledgements.} The author would like to sincerely thank Charles Torossian for his inspiring ideas and guidance during the preparation of his thesis. He would also like to thank the EU RTN Liegrits and its coordinator Fred Van Oystaeyen for financial support.

\end{document}